\input amstex
\documentstyle{amsppt}
\vsize=8in \hsize 6.6 in 
\loadbold \topmatter
\title The Closed Orbit Controllability Criterium \endtitle
\author Valeri Marenitch \endauthor

\address Kalmar, Sweden\endaddress \email valery.marenich\@ gmail.com \endemail
 \subjclass 93B05, 93B29, 20M20 \endsubjclass
\date{}\enddate
\abstract We prove that every closed "general" trajectory of the control system $\Sigma_M$ has an open neighborhood on which $\Sigma_M$ is
controllable if 1) this orbit contains some point where the Lie algebra rank condition ($LARC$) is satisfied, and 2) the set of control vectors is
"involved" at $q$. In particular, for the control systems $\Sigma_M$ on the compact connected manifold $M^n$ with an open control set this gives
the following "Closed Orbit Controllability Criterium": The dynamical system $\Sigma_M$ of the considered type is controllable on $M^n$ if and
only if for an arbitrary point $q$ of $M^n$ there exists a closed trajectory of the control system going through this point. We also present
examples which show that our conditions are necessary.
\endabstract
\dedicatory {Preliminary version. Comments are welcome.}\enddedicatory
\endtopmatter

\document

\head Introduction and results \endhead

In this note we study the chronological map ${\Cal G}(t_1,...,t_N): R^N\to M^n$ of the control dynamical system
$$
\dot q(t)=V_u(q(t)) \qquad \qquad \qquad \qquad (\Sigma_M)
$$
on a complete manifold $M^n$. We prove that it is an open map in the open domain in $R^N$ for sufficiently big $N$, if the system satisfies the
well-known Lie Algebra Rank condition ($LARC$). This implies the controllability of the system in some neighborhood of any, so called "general"
closed orbit.

To ensure that an arbitrary trajectory of the system can be approximated by some general trajectory we use the "involvement" condition which,
essentially, means that any direction $V_u$ lies in some positive cone generated by some set of $V_{u_i}$ satisfying $LARC$. More precisely, we
prove that if the set of control vectors $V_u$ is open in the linear subspace $L$ which it generates (the "open" condition), and its convex hull
coincides with $L$ (the "ample" condition), then $\Sigma_M$ is "involved", see definitions below.

For the control dynamical systems $\Sigma_M$ on a compact connected manifold $M^n$ which

\medskip

1) satisfy the Lie algebra rank condition everywhere and

2) have "involved" control vectors set;

\medskip

we prove the following "Closed Orbit Controllability Criterium": such system is controllable if and only if
through every point of the manifold goes some closed orbit of $\Sigma_M$, see Theorem~A below.

\medskip

One of the applications of the Theorem~A is to bilinear systems in the Euclidean space. These are the systems of
linear equations
$$
\dot x(t)= A(u)x(t) = (A + u^1 B_1 + ... +u^d B_d)x(t), \qquad \qquad \qquad \qquad (\Sigma)
$$
where the right-hand side - the linear operator $A(u)$ - is a linear combination of some constant linear operators $A$ and $B_k$ which do not
depend on the control parameter $u\in R^d$. For such systems the "Closed Orbit Controllability Criterium" provides many new conditions for
controllability, both necessary and sufficient; see our forthcoming paper [CM]. Here we mention only the simplest one for the systems in
three-dimensional Euclidean space $R^3$, see Theorem~B below: if for two control parameters $u$ and $v$ the right-hand side linear operators have
complex eigenvalues $\lambda_C(u)$ and $\lambda_C(v)$ correspondingly, then the system $\Sigma$ satisfying $LARC$ is controllable if for the real
eigenvalues $\lambda_R(u)$ and $\lambda_R(v)$ of $A(u)$ and $A(v)$ it holds
$$
( \lambda_R(u) - Re(\lambda_C(u)) ) \quad ( \lambda_R(v) - Re(\lambda_C(v)) ) < 0.
$$

The results of this paper were obtained during a visit to the Universidad Michoacana (Morelia, Mexico). The author
sincerely thanks Prof A.~Choque for the organization of this visit, cooperation and hospitality.

\medskip

\head 1. Twisted LARC \endhead

To introduce notations we remind some basics of the theory of control systems, see [AS].

Let $M^n$ be a complete n-dimensional manifold without boundary. The control system on $M$ is given by the family
of dynamical systems
$$
\dot q(t)=V_u(q(t)), \quad q\in M, u\in U; \qquad (\Sigma_M),  \tag 1
$$
where the right-hand side is a smooth vector field $V_u$ depending on the control parameter $u$ from some set $U$
of controls (usually the domain in some Euclidean space: $U\subset R^d$). The solution of $\Sigma_M$ is the
trajectory $q(t)$ such that
$$
\dot q(t)=V_{u(t)}(q(t)), \quad q\in M, u(t): R\to U \tag 2
$$
for some admissible control function $u(t)$ (usually, the locally integrable functions $u(t):R\to U\subset R^d$).

\medskip

\subhead 1.1. Semi-group ${\Cal G}^+$ of $\Sigma_M$ \endsubhead

For a given $u\in U$ denote by $P^t_u$ the one-parameter group of diffeomorphisms of $M$  generated by the vector
field $V_u$ on $M$. The standard notations are:
$$
P^t_u = P^t_{V_u} = exp(\int\limits_0^t V_u d\tau)=e^{tV_u},  \tag 3
$$
where $t$ might be positive, negative or zero (and then $P^0_u$ is the identity map $Id$). Compositions of such
diffeomorphisms for different vector fields $V_{u_i}$, called "chronological products",
$$
{\Cal G}(t_1,t_2,...,t_N)=P^{t_1}_{u_1}\circ ... \circ P^{t_N}_{u_N} \tag 4
$$
generate the Group ${\Cal G}={\Cal G}(\Sigma_M)$ of the system, which is the (Frechet-) subgroup of the complete group $Diff(M)$ of all
diffeomorphisms of $M$. If we restrict ourself only to nonnegative parameters $t_i$ we have the (positive) semigroup ${\Cal G}^+$ of the system.
If we assume that every control function $u(t)$ can be approximated by a piece-wise constant function $u_N(t)$ which equals $u_i$ on the interval
$(T_{i-1},T_i)$ where $T_i=t_N+...+t_{i}$; and that every trajectory of the control system $\Sigma_M$ can be approximated by trajectories with the
piece-wise controls; then  the set ${\Cal O}^+(q)$ of points reachable from $q$ by trajectories of $\Sigma_M$ equals the orbit ${\Cal G}^+(q)$ of
$q$ under the action of the semigroup ${\Cal G}^+$.

The system $\Sigma_M$ is controllable if for every two points $q,p$ there exists some trajectory $q(t,u(t)), 0\leq
t\leq T$ of the system $\Sigma_M$ starting at $q=q(0,u(0))$ and ending at $p=q(T,u(T))$. In other words, the
system is controllable if the action of the semi-group ${\Cal G}^+$ is transitive. Since it is easier to handle
the group-actions than actions of semi-groups, the question of whether $\Sigma_M$ is controllable or not is,
usually, divided into two sub-questions:

1) whether the action of the group ${\Cal G}$ is transitive or not,

2) whether the orbits of the semi-group ${\Cal G}^+$ actually coincide with the orbits of the bigger ${\Cal G}$.

Although in general the first question is not easy, its local version though has the following well-known answer.
We may call the system $\Sigma_M$ {\bf locally transitive at $p$} if there exists some open
$\epsilon(q)$-neighborhood $B_{lt}(q)$ such that every point $p$ from this neighborhood if reachable from $q$ by
some trajectory $q(t,u(t)), 0\leq t\leq t(q,\epsilon(q))$ of the $\Cal G$-action, where $t(q,\epsilon(q))\to 0$ as
$\epsilon(q)\to 0$. Then the local transitivity in $q$ is assured by the well-known Lie algebra rank condition
$(LARC)$\footnote{such systems are also called completely nonholonomic or bracket-generating}:

1') {\bf $LARC(q)$}  The Lie-algebra ${\Cal LA}$ generated by some number of vector fields $V_{u_i}, i=1,...,d'$
is such that its evaluation at the point $q$ coincide with the whole tangent space $T_qM$ to $M$ at $q$.

Indeed, note that in general the one-parameter diffeomorphisms in the "chronological product" do not commute.
Thus, it is possible to move the point by the action of ${\Cal G}$ not only along vector fields $V_{u_i}$ for
different $i$, but also along their commutators, see the (11) below; and so on. Therefore, the orbit ${\Cal G}(q)$
- by The Orbit Theorem, see [AS] - is some immersed submanifold in $M^n$ with the space of tangent vector fields
closed under the Lie bracket. Then from the $LARC$ at $q$ we see that the tangent space of the orbit ${\Cal G}(q)$
equals the tangent space to the whole $M^n$, and that the orbit ${\Cal G}(q)$ contains some open neighborhood of
$q$. The last means that the action of $\Cal G$ is locally transitive.\footnote{Of course, the actual choice of
"generators" $V_{u_i}, i=1,...,m$ for systems $\Sigma_M$ which satisfy $LARC$ is not canonical, but depends on the
particular form of $V_u$. For instance, when the vector field $V_u$ depends smoothly on $u$ we have
$$
V_u=V_{u_0} + dV_{u_0}(u-u_0) + o(\|u-u_0\|) = A + (u-u_0)^i B_i + o(\|u-u_0\|)
$$
for $A(q)=V_{u_0}(q)$ and $B_i(q)= \partial V_u(p)/\partial u^i$ at every $u_0$. Then for the smooth system
satisfying  $LARC$ we may find $u_0\in U, U\subset R^d$ such that the Lie-algebra generated by the vector fields
$A$ and $B_i$ for $u_i=u_0+\delta E_i$ for all sufficiently small $\delta$, where $E_i$ is the unit basis vector
in $R^d$ with number $i$; coincide with the Lie-algebra generated by vector fields $A$ and $B_i$ at $u_0$. Then,
without loss of generality we may assume, that the "generators" above are the vector fields $A,B_i$. We may
reformulate this by saying that under the $LARC$ the tangent space $T_qM^n$ to $M^n$ at the point $q$ is generated
by vectors $A(q)$, $B_i(q)$ and their commutators: say $C_{ij}(q)=[B_i,B_j](q)$, and so on. This happens in the
particular case of so called bilinear systems in Euclidean spaces when $M^n=R^n$:
$$
\dot q(t)=(A + u^1 B_1 + ... +u^d B_d) q(t),
$$
which we consider in [CM].}

In general, the point $p$ close to $q$ could be reached by "long" trajectories issuing from $q$ so that the
$LARC(q)$ at $q$ is not necessary for controllability or transitivity of the $\Cal G$-action, see the example 5.3
in [AS]. Actually, slightly generalizing this example it is easy to construct controllable dynamical systems such
that the $LARC$ does not hold at any point at all, i.e., such that ${\Cal LA}(q)\not= T_qM^n$ everywhere: say, for
arbitrary $n>2$ we may always have $dim({\Cal LA}(q))\leq 2$.

\medskip

\proclaim{Example~1} For an arbitrary $n>2$ there exists a controllable dynamical system $\Sigma_T$ on
$n$-dimensional torus $T^n$ such that in all points $dim{LA}(q)\leq 2$.
\endproclaim

\demo{Construction}

Indeed, take the product $\Pi^n \subset R^n$ of the closed interval $[0,n-1]$ with $(n-1)$-dimensional unit cube
$[0,1]^{n-1}$, and introduce coordinates $\overrightarrow{x}=\{x,y^2,...,y^n\}$ in $\Pi$ correspondingly. Choose
some smooth nonnegative function $\omega(r)$ with support inside $[1/3,2/3]$ and an integral bigger than $1$, and
define the vector field $V(x,y^2,...,y^n,u)$ in the parallelepiped $\Pi$ as follows: for $k\leq x\leq k+1$ where
$k=0,1,...,n-2$ is integer take
$$
V(x,y^2,...,y^n,u) = {{\partial}\over{\partial x}} + u \omega(x-k) {{\partial}\over{\partial {y^{k+2}}}}, \tag
exm~1
$$
where the control $u$ belongs to an open interval ${\Cal \Omega}=(-1, 1)$. Since $V$ is invariant under parallel
translations of $R^n$ of the form:
$$
V(x,y^2,...y^k,...,y^n,u)=V(x,y^2,...y^k+1,...,y^n,u)
$$
and
$$
V(x,y^2,...y^k,...,y^n,u)=V(x+d,y^2,...y^k,...,y^n,u);
$$
- in fact, it depends only on $x$ - we see that after identification of opposite sides of the parallelepiped $\Pi$
we obtain the $n$-dimensional torus $T^n$ with a smooth vector field which we denote again by
$V(\overrightarrow{x},u)$. All trajectories of $V(\overrightarrow{x},u)$ are curves normal to the sub-torus
$T^{n-1}$ defined by $x=0$ with $\{y^2,...,y^n\}$ being local coordinates in $T^{n-1}$. From the definition above
we see that changing the control $u(t)$ on the interval $t\in [k-1,k], k=0,...,n-2$ leads to the change of the
$y^{k+2}$-coordinate, leaving all other coordinates intact. Thus, the point $p$ with arbitrary coordinates
$\{y^2,...,y^n\}$ can be reached from the zero-point $q$  (possibly, after one additional "rotation" with zero
control along $x$-coordinate parallel of the torus $T^n$) with the help of the following piece-wise constant
control function $u(t)$:
$$
u(t) = y^k (\int \omega(r) dr)^{-1} \tag exm~2
$$
on the interval $k-1\leq t\leq k$. We see that ${\Cal O}^+(p)=T^n$, and the control system
$$
{{d}\over{dt}} \overrightarrow{x} = V(\overrightarrow{x},u) \qquad \qquad \qquad \qquad (\Sigma_T) \tag exm~3
$$
is controllable, while the dimension of ${\Cal LA}$ in each point is 1 or 2.
\enddemo

\medskip

Despite the example above, the controllability implies some generalized $iLARC$ ("integral"). We address this question later, while here we note
only the following: if the action of ${\Cal G}$ is transitive then the orbit of every point under this action is the whole manifold $M^n$, i.e.,
is an open subset in particular. We know that every orbit is an immersed manifold. Therefore, since $M^n$ is a second category set; it has to have
the same dimension $n$ as $M^n$. Then compositions (4) ("chronological products") with some fixed $N$ should provide us with an open map
(surjection) from the space of $(t_1,...,t_N)$ to some neighborhood of $q$ in $M^n$. The standard arguments then show that controllability leads
to some kind of generalized $LARC$ along all closed trajectories through $q$ compared with the help of corresponding Poincare maps.

\medskip

In this paper we restrict ourself to the case when the $LARC$ holds at the point $q$ we consider:
$$
{\Cal LA}(q)=T_q M^n, \tag 5
$$
and then prove that every "general" closed orbit through $q$ belongs to the orbit ${\Cal G}^+(q)$ together with some its open neighborhood
$B^+_q$, and further, has some open neighborhood $B_q$ on which the system $\Sigma_M$ is controllable.

\medskip

\subhead 1.2. Vector fields as operators on $C^\infty(M)$ (after [AS])\endsubhead

We know that the manifold $M^n$ may be identified with homomorphisms $\phi: C^{\infty}(M)\to R$ of the algebra of
smooth function on $M$; i.e., for each algebra homomorphism $\phi$ there exists some point $q\in M$ such that
$\phi$ coincides with the evaluation homomorphism $\hat{q}(a)=a(q)$ at this point, see [AS].\footnote{Note that in
the proof of this (Gelfand-Neimark type) theorem in the Lemma~A1 in [AS] the set $M\backslash K$ can be empty.}

In the same way it is useful to represent diffeomorphisms $P:M\to M$ as automorphisms of the algebra $C^\infty(M)$
of smooth functions on $M$ acting by change of variables: for any given function $a$ on $M$ the action $\hat{P}$
on $a$ is given by $\hat P(a)=a\circ P$. Indeed, every algebra homomorphism $A:C^\infty(M)\to C^\infty(M)$ is some
$\hat P$ defined by $P(q)=q_1$ where $\hat q_1=\hat q\circ A$.

Then the vector field $V$ is given by the derivation of $C^\infty(M)$, i.e., the linear operator $\hat V:
C^\infty(M)\to C^\infty(M)$ satisfying the Leibniz rule $\hat V(ab)=\hat V(a)b+a\hat V(b)$. If the vector field
$V$ generates the flow $P^t$ then from the definition it follows
$$
{{d}\over{dt}}\hat P^t =\hat V\circ P^t = P^t\circ \hat V. \tag 6
$$
If $P_*$ denote the differential of $P$ then the vector field $P_*(V)$ equals the derivation $P^{-1}\circ V\circ
P$. This justifies the notation $P_*=Ad(P)$, while the calculation
$$
{{d}\over{dt}} (Ad P^t)W_{| t=0} = {{d}\over{dt}} (P^t\circ W \circ P^{-t})_{| t=0} =[V,W] \tag 7
$$
shows that the Lie bracket of vector fields $V$ and $W$ is given by the derivative of $Ad(P^t)$ at $t=0$ which is
denoted by $(adV)W=[V,W]$. Since the Lie bracket is preserved under the $P_*$ we conclude
$AdP^t[X,Y]=[AdP^tX,AdP^tY]$ which after differentiation at $t=0$ gives the Jacobi identity
$$
(adV)[X,Y]]=[(adV)X,Y] + [X,(adV)Y] = [V,[X,Y]]=[[V,X],Y] + [X,[V,Y]].
$$
From (6) we may conclude the following asymptotic form for $P^t$, also denoted as
$$
exp({tV}) = \sum\limits_{n=0}^{\infty} {{t^n}\over{n!}} V^n, \tag 8
$$
where $V^n=V\circ ... \circ V$ (one variable conjugate Taylor formula). If
$$
S_m(t)= Id + \sum\limits_{n=1}^{m-1} {\int} ... \int\limits_{\Delta_n(t)} V_{\tau_n}\circ ... \circ V_{\tau_1}
d\tau_n ... d\tau_1 \tag 9
$$
is the partial sum, then for any $a\in C^\infty (M), s\geq 0$ and compact $K\subset M$ it holds
$$
\|(exp\int\limits_0^t V\tau d\tau - S_m(t))a \|_{s,K} \leq O(t^m), \tag 10
$$
see (2.13) from [AS]. Below we assume that $m$ is big enough, say $m>2^n$; which is sufficient for estimates
below. Using (9-10) for diffeomorphisms generated by different vector fields we deduce the asymptotic
representation for an arbitrary chronological product (the conjugate Taylor formula for many variables), or
approximation formula for the action of the group ${\Cal G}$. For instance, for the curve $q(t)=q\circ P_v^t\circ
P^t_w\circ P_v^{-t}\circ P^{-t}_w$ it gives the following expansion:
$$
q(t)= q\circ (Id + t^2 (V\circ W - W\circ V) + o(t^2)), \tag 11
$$
where $V,W$ are generators of the one-parameter groups of diffeomorphisms $P_v^t$ and $P^t_w$. Generally, for the
chronological product ${\Cal G}(t_1,...,t_N)=P^{t_1}_{u_1}\circ ... \circ P^{t_N}_{u_N}$ the same arguments give:
$$
q(t_1,...,t_N)=q\circ {\Cal G}(t_1,...,t_N) = q\circ (Id +  \sum\limits_{n=1}^{m-1} {{t_1^n}\over{n!}} V_{u_1}^n) \circ ... \circ (Id +
\sum\limits_{n=1}^{m-1} {{t_N^n}\over{n!}} V_{u_N}^n) + O(T^m), \tag 12
$$
where by $T$ we denote $|t_1|+...+|t_N|$. We may further expand (open all brackets) the above formula (12) to
obtain the Taylor type approximation:
$$
q(t_1,...,t_N)= q\circ (Id + \sum_{k=1}^{k=m-1} \sum_{\alpha^k}
{{t_1^{\alpha_1^k}...t_N^{\alpha_N^k}}\over{\alpha_1! ... \alpha_N!}} V_{u_1}^{\alpha_1}\circ ...\circ
V_{u_n}^{\alpha_N}) + O(T^m), \tag 13
$$
where $\alpha^k=(\alpha_1^k,...,\alpha_N^k)$ is $N$-multi-index of order $k$; i.e., $\|\alpha^k\|=\alpha_1^k+...+\alpha_N^k=k$; while
$V_{u_i}^{\alpha_i}=V_{u_i}\circ ...V_{u_i})$ denotes the composition of $\alpha_i$ derivations $V_{u_i}$. For the notational convenience, we
consider the vector fields $V_{u_i}$ for $i=1,...,N$ to be different, also they may be (and usually are) taken from some finite set $V_i,
i=1,...,d'$ of vector fields (say, satisfying the $LARC$). It is also important that, despite the fact that the asymptotic representations of
chronological products above in general may be not convergent for some collections of vector fields $V_u$ and on some functions from
$C^\infty(M^n)$, in our case - when the set of generators $V_{u_i}$ is finite and we consider the action on coordinate functions - representation
above are convergent (and are, actually, nothing more than the Taylor formulas).\footnote{In fact, we may use these finite approximations for the
chronological products as their substitutes to obtain all our results below and the Theorem~A in particular. The estimate (10) then guarantee that
the claim of the Theorem~A is valid for the asymptotic limit under our conditions - openness of $V(q)$ - on the set of control vectors.}

\medskip

\subhead 1.3. Rank of the chronological map \endsubhead

Due to (6) for the given chronological product $q(t_1,...,t_N)$ its partial derivatives are:
$$
{{\partial}\over{\partial t_i}}q(t_1,...,t_N) = {{\partial}\over{\partial t_i}}P^{t_1}_{u_1}\circ ... \circ P^{t_N}_{u_N} = P^{t_1}_{u_1}\circ
...\circ V_{u_i}\circ P^{t_i}_{u_i}\circ ... \circ P^{t_N}_{u_N}, \tag 14
$$
which with the help of (9,10) leads to the following asymptotic representations:
$$
{{\partial}\over{\partial t_i}}q(t_1,...,t_N) = {{\partial}\over{\partial t_i}} q\circ (Id +
\sum\limits_{n=1}^{m-1} {{t_1^n}\over{n!}} V_{u_1}^n) ...\circ V_{u_i}\circ (Id +  \sum\limits_{n=1}^{m-1}
{{t_i^n}\over{n!}} V_{u_i}^n) ... \circ (Id + \sum\limits_{n=1}^{m-1} {{t_N^n}\over{n!}} V_{u_N}^n) + O(T^m), \tag
15
$$
or, comparing with (13)  we get:
$$
{{\partial}\over{\partial t_i}}q(t_1,...,t_N) = q\circ \sum_{k=1}^{k=m-1} \sum_{\alpha^k} (
{{\partial}\over{\partial t_i}} {{t_1^{\alpha_1^k}...t_N^{\alpha_N^k}}\over{\alpha_1! ... \alpha_N!}} )
V_{u_1}^{\alpha_1}\circ ...\circ V_{u_N}^{\alpha_N} + O(T^{m-1}). \tag 16
$$

Again, we remind that the chronological map ${\Cal G}:(t_1,...,t_N)\to q(t_1,...,t_N) \in M^n$ we understand here
as the map from the Euclidean space $R^N$ to algebra homomorphisms of $C^\infty(M)$, i.e., each point of $M$ is
understood as the evaluation homomorphism. Then partial derivatives of $q(t)$ above are derivations of
$C^\infty(M^n)$ corresponding to the vector fields. In a "usual sense" - as sections of the tangent bundle of
$M^n$ - these partial derivatives are
$$
W_i(q(t_1,...,t_N))=d{\Cal G}(t_1,...,t_N)({{\partial}\over{\partial t_i}}), \tag 17
$$
where here the chronological map ${\Cal G}$ is the smooth map from the Euclidean space $R^N$ to the differential
manifold $M^n$.

Now we introduce some local coordinates $(x^1,...,x^n)$ in some small neighborhood which contains both points
$q'=q(t_1,...,t_N)$ and $q$ - assuming of course that they are close enough; i.e., that $|T|$ is small. In this
neighborhood we have coordinate functions $x^j: M^n\to R$, which we may assume to be continued outside the
coordinate neighborhood in an arbitrary way; i.e., it holds: $x^j\in C^\infty(M^n)$. The corresponding coordinate
vector fields $E_i$ are such that
$$
E_i(x^j)\equiv \delta_i^j \tag 18
$$
in this coordinate neighborhood. Then for an arbitrary differential operator $D^\alpha$ of order $|\alpha|$
strictly bigger than one - not having order zero and one terms - it holds
$$
D^\alpha(x^j) \equiv 0. \tag 19
$$
Each vector $W$ tangent to $M^n$ at the point $q'=q(t_1,...,t_N)$ is defined by its action on coordinate functions
as follows:
$$
W=\sum\limits_{i=1}^{n} (W(x^i)) E_i. \tag 20
$$

Rewrite (16) separating differential operators by their orders:
$$
{{\partial}\over{\partial t_i}}q(t_1,...,t_N) = q\circ (\sum\limits_{k=1}^{m-1} d_{i\alpha}^k(t_1,...,t_N)
D_{\alpha}^k(V_{u_1},...,V_{u_N})) + O(T^{m-1}), \tag 21
$$
where $d_{i\alpha}^k$ are some polynomials of $t_i$,  and $D_{\alpha}^k$ - some finite number of differential
operators of order $k$, i.e., such that in our coordinate system
$$
D^k_{\alpha}= c^k_{\alpha,s}(x_1,...,x_n){{\partial^k}\over{\partial^{{s_r}} x_{r}}}, \tag 22
$$
where $c^k_{\alpha,s}(x_1,...,x_n)$ - some smooth functions, and ${s_r}=({s_1},...,{s_n})$ some multi-index of order $k$ - ${s_1}+ ... + {s_n}=k$.
Due to (19) it holds
$$
{{\partial}\over{\partial t_i}}q(t_1,...,t_N)(x^j)= d^1_{i\alpha}(t_1,...,t_N) D^1_{\alpha}(x_j) + O(T^{m-1}).
\tag 23
$$
In other words, we see that the differential of the chronological map can be recovered from its action on linear
coordinate functions, while this action is determined by its "first order operator part". More precisely, the
differential $d{\Cal G}(t_1,...,t_N)$ of the chronological map at the point $q$ on the vector
${{\partial}\over{\partial t_i}}$ equals the derivation - the vector field we denote $F_i(t_1,...,t_N)$, which
action on the coordinate function $x^j$ is given by:
$$
F_i(t_1,...,t_N)(x^j) = (d{\Cal G}(t_1,...,t_N) ({{\partial}\over{\partial t_i}}))^j = {{\partial}\over{\partial
t_i}}q(t_1,...,t_N)(x^j)= d^1_{i\alpha}(t_1,...,t_N) D^1_{\alpha}(x_j) + O(T^{m-1}). \tag 24
$$
Examine $F_i$ more closely. From (16) we conclude:
$$
F_i(t_1,...,t_N))= \sum_{k=1}^{k=m-1} \sum_{\alpha^k} {{\partial}\over{\partial t_i}}
({{t_1^{\alpha_1^k}...t_N^{\alpha_N^k}}\over{\alpha_1! ... \alpha_N!}}) D^1( V_{u_1}^{\alpha_1}\circ ...\circ
V_{u_N}^{\alpha_N}) + O(T^{m-1}), \tag 25
$$
where by $D^1_\alpha=D^1(V_{u_1}^{\alpha_1}\circ ...\circ V_{u_n}^{\alpha_N})$ we denote the first differential
operator part of the compositions of derivations $V_{u_i}$; which by direct calculations equals
$$
D^1_\alpha = D^1(V_{u_1}^{\alpha_1}\circ ...\circ V_{u_N}^{\alpha_N}) = (v_1^{s^1_1} v_1^{s^1_2} ...
v_1^{s^1_{\alpha_1}} {{\partial^{\alpha_1}}\over{\partial x^{s_1}...\partial x^{s_{\alpha_1}}}} (v_2^{s^2_1} ...
(v_N^{s^N_j} )...) {{\partial}\over{\partial x^j}}, \tag 26
$$
where $v^s_i$ are coordinates of the vector field $V_{u_i}$:
$$
V_{u_i}=v^s_{i} {{\partial}\over{\partial x^s}}. \tag 27
$$
Comparing (24) with (25) we see that polynomial coefficients of first order differential operators $D^1_\alpha$
(which (operators) do not depend on $(t_1,...,t_N)$) are the following monomials:
$$
d^1_{i\alpha}(t_1,...,t_N)={{\partial}\over{\partial t_i}} ({{t_1^{\alpha_1^k}...t_N^{\alpha_N^k}}\over{\alpha_1!
... \alpha_N!}}). \tag 28
$$
For what follows it is important that such monomials are linearly independent; i.e., if some linear combination
$\sum\limits_\alpha \Lambda_\alpha d^1_{i\alpha}(t_1,...,t_N)$ is identically zero, then all coefficients
$\Lambda_\alpha$ are zeroes.\footnote{To prove that in such combination $\Lambda_\alpha=0$ take the partial
derivative $\partial^\alpha/\partial t_1^{\alpha_1}...$ of this linear combination at zero.} If some
$D^1_\alpha(x^j)\not= 0$ in this linear combination, then the zero set of the action on $x^j$ of the first order
differential part of the differential of the chronological map at the point $q$ on the vector
${{\partial}\over{\partial t_i}}$ vanishes on the zero set given by the polynomial equation:
$$
\Omega^j_i = \{(t_1,...,t_N) | \sum\limits_\alpha d^1_{i\alpha}(t_1,...,t_N) D^1_\alpha(x^j) =0 \}. \tag 29
$$

As the corollary we conclude the following.

\proclaim{Lemma~1} The zero set $\Omega^j_i\subset R^N$ given by the polynomial equation is closed and has an
empty interior.
\endproclaim

\medskip

Next we note that due to our asymptotic expansions (23), (24) for sufficiently small $|T|$ the differential of the
chronological product ${\Cal G}(t_1,...,t_N): q\to q(t_1,...,t_N)$ is not a surjection on $R^N$ at some point
$q'(t')$ for $t'=(t'_1,...,t'_N)$ if for some system of coordinates of the type we consider, the action of all
$D^1_\alpha$ on some particular function $x^j$ vanish at zero; i.e., such that
$$
q\circ D^1_\alpha (x^j) = 0. \tag 30
$$
Indeed, if the differential of ${\Cal G}$ is not surjection we may construct (using the Orbit theorem) the coordinate system $x^j, 1,...,n$ in
some neighborhood of $q$ such that the first $n'<n$ coordinate vectors $E_i, i=1,...,n"$ at the point $q'=q(t'_1,...,t'_N)$ generate the image of
the differential of ${\Cal G}$ at the point $q$ - some $n"$-dimensional subspace of $T_q'M^n$ with $n"<n$. Then the derivatives of all $x^j, j>n"$
in all direction of this subspace vanish, i.e., (30) holds. If otherwise; i.e., we are given some chronological product ${\Cal G}(t_1,...,t_N):
R^N\to M^n$ such that (30) is not satisfied for an arbitrary coordinates $x^i, i=1,...,n$ then this map is a surjection on the open subset
$R^N\backslash \Omega$ which complement $\Omega$ has an empty interior, for $|T|< \epsilon$ for some small $\epsilon$. Since from (23) we see that
$$
q\circ D^1_{\alpha}(x_j)  = {{\partial^k}\over{\partial t_1^{\alpha_1}... \partial t_i^{\alpha_i} ...\partial
t_N^{\alpha_N}}} {\Cal G}(t_1,...,t_N)_{|t_1=...=t_N=0}(x^j),  \tag 31
$$
combining our arguments above we come to the following statement.

\proclaim{Lemma~2 ("PD-LARC")} The chronological map
$$
{\Cal G}_{(u_1,...,u_N)}(t_1,...,t_N)=P^{t_1}_{u_1}\circ ... \circ P^{t_N}_{u_N} \tag 32
$$
is a surjection for $|T|=|t_1|+...+|t_N|<\epsilon$, where $\epsilon$ is sufficiently small, and for
$t=(t_1,...,t_N)\in R^N \backslash \Omega$ where $\Omega$ closed set with an empty interior, given by some
polynomial equation; if there exists $n$ linearly independent partial derivatives
$$
D_j = D^1 ({{\partial^\alpha(j)}\over{\partial t_1^{\alpha(j)_1^k} ...\partial t_N^{\alpha(j)_N^k}}} {\Cal
G}(t_1,...,t_N)_{|t_1=...=t_N=0})(x^j), \tag 33
$$
where $D^1$ denotes the first order differential part, and $\alpha(j)$ is some multi-index $(\alpha(j)_1,...,\alpha(j)_N)$.
\endproclaim

\medskip

Note, that the Lemma~2 condition is not stronger than the $LARC$. If the $LARC$ is not satisfied then, using the
orbit theorem, we may always construct the coordinate system in $M^n$, such that the orbit of $q$ is contained in
some coordinate subspace given by, say $x^n=0$. Then all partial derivatives of ${\Cal G}$ will have the first
order differential parts with the vanishing action on $x^n$.

Now we are ready to construct the chronological maps as in the Lemma~2 above. Since any exponential $P_u^t$ is the
identity map when $t=0$, the partial derivative of (32) in (33) equals the same partial derivative of the
chronological product only of those $P^{t_i}_{u_i}$ with $\alpha(j)_i \not= 0$. Therefore, our chronological
product satisfying conditions of the Lemma~2 may be constructed as the superposition of the chronological products
$$
{\Cal G}^j(t^j_1,...,t^j_{k(j)})= P^{t^j_1}_{u^j_1} \circ ... \circ P^{t^j_{k(j)}}_{u^j_{k(j)}}, \qquad t^j=(t^j_1,...,t^j_k(j)) \in R^{k(j)},
\tag 34
$$
which partial derivative has the first order differential part equals some $E_j$ such that $E_j,j=1,...,n$ are
linearly independent. Now the $LARC$ guarantees the existence of such ${\Cal G}^j$.

Indeed, if $\Sigma_M$ satisfies the $LARC$ we may assume that the basis $E_j,j=1,..,n$ consists of some number of
vectors $V_{u_i}, i=1,..., d'$ - which we denote by $W_i, i=1,..,d'$ and their mutual commutators $W_i,
i=d'+1,...,n$, where
$$
W_i=[W_{\beta^i_k}, [W_{\beta^i_{k-1}}, [... ]]]; i=d'+1,...,n \tag 35
$$
for some $k$-multi-index $\beta^i=(\beta^i_k,...,\beta^i_1)$, where the order of the commutator $k$ in turn
depends on the index $i$ of the vector $W_i$: i.e., $k=k(i)$. Easy to see that orders of commutators are bounded
by $n$ - the dimension of $M^n$.\footnote{Hence $m>2^n$ is sufficient.}

By definition (35) the derivation $W_i$ is the linear combination of the derivations which are the first order
differential parts
$$
W_i= \sum_\sigma (-1)^{\sigma} D^1( W_{\beta^i_\sigma(k)}\circ W_{\beta^i_{\sigma(k-1)}} ... ) \tag 36
$$
of differential operators of the form
$$
W_{\beta^i_\sigma(k)} \circ W_{\beta^i_{\sigma(k-1)}} \circ ... \circ W_{\beta^i_{\sigma(1)}}, \tag 37
$$
where $\sigma$ some permutations of indexes $\beta^i_j$. An arbitrary derivation in the linear combination (36) is
the partial derivative of the chronological product of the type:
$$
W_{\beta^i_\sigma(k)} \circ W_{\beta^i_{\sigma(k-1)}} \circ ... \circ W_{\beta^i_{\sigma(1)}} = {{
\partial^k}\over {\partial t_k \partial t_1}}={\Cal G}^i_{\sigma}(t_1,...,t_k)=
P^{t_1}_{\beta^i_\sigma(k)} \circ P^{t_{k-1}}_{\beta^i_{\sigma(k-1)}} \circ ... \circ
P^{t_1}_{\beta^i_{\sigma(1)}}. \tag 38
$$
Taking composition of such products for all $i=1,...,n$ with different arguments $t_s$ we obtain (very long!)
chronological products whose partial derivatives realize all monomials in (36). We call such chronological product
{\bf the general at $q$  relative to the set of generators $\{V_1,...,V_{l'}\}$}. They satisfy the Lemma~2, which
may be reformulated as follows.

\proclaim{Lemma~3} If $\Sigma_M$ satisfies the $LARC$ at $q$ then there exist the general chronological product
${\Cal G}:R^N\to M^n$ at $q$, such that its differential is a surjection on $R^N\backslash \Omega$ for some closed
$\Omega\subset R^N$ with an empty interior which is given by some polynomial equation.
\endproclaim

\medskip

\proclaim{Definition} The trajectory $q(t), 0\leq t \leq \epsilon$ of the control system $\Sigma_M$ is called {\bf
the general relative to the set of generators $\{V_1,...,V_{l'}\}$} if it contains some interval which goes from
$q=q(0)$ to some point $q(\epsilon)={\Cal G}(t_1,...,t_N)(q)$ with all $t_i$ positive and such that
$t=(t_1,...,t_N)\in R^N\backslash \Omega$. \endproclaim

\medskip

If $q(t), 0\leq t \leq \epsilon$ is the general trajectory relative to some set of generators of $\Sigma_M$ then
by definition some general chronological product ${\Cal G}$ provides the surjection from some neighborhood of zero
in $R^N$ over some neighborhood $B(\epsilon)$ of $q(\epsilon)$, which means that all points from this neighborhood
are reachable from $q$, or $B(\epsilon)\subset {\Cal O}^+(q)$. If further, this interval $q(t), 0\leq t \leq
\epsilon$ is a part of some closed trajectory $q(t), 0\leq t \leq t_q$ going through the point $q=q(t_q)$ again
under some control $u(t): [0,t_q] \to {\Cal U}$ then by the standard continuous dependence arguments the image of
this open neighborhood $B(\epsilon)$ under the flow $F(t,u(t))$ is the open neighborhood of $q(t), \epsilon < t$
also reachable from $q$. The union of all these neighborhoods gives us the open neighborhood $B^+_q$ of the closed
trajectory $q(t), 0 \leq t\leq t_q$, which by construction belongs to the positive orbit of $q$: $B^+_q \subset
{\Cal O}^+(q)$.

We formulate the obtained results as follows.

\proclaim{Theorem~1} If the control dynamical system $\Sigma_M$ satisfies the $LARC$ at the point $q$; i.e., there
exists the complete system - the set of generators -  $E_j,j=1,..,n', n'\geq n$ ($E_j$ generate $T_qM^n$) of
$T_qM^n$ consisting of some number of vectors $V_{u_i}, i=1,..., d'$, denoted by $W_i, i=1,..,d'$; and their
mutual commutators $W_i, i=d'+1,...,n'$
$$
W_i=[W_{\beta^i_k}, [W_{\beta^i_{k-1}}, [... ]]]; i=d'+1,...,n',
$$
then the set of points $q(t_1,...,t_N)$ (where all $t_i$ are positive) reachable from $q$ by general trajectories relative to this set of
generators is an open subset of $M^n$. If there exists a closed trajectory $q(t), 0\leq t \leq t_q$ going through the point $q$ then some open
neighborhood $B^+_q$ of this trajectory is reachable from $q$:
$$
q' \in {\Cal O}^+(q) \quad \text{ for all } \quad q' \in B^+_q.
$$

\endproclaim

Now consider the negative orbit ${\Cal O}^-(q)$ of $q$ - the chronological products ${\Cal G}(t_1,...,t_N)$ where all $t_i$ are non-positive. By
the definition these are the points from which we may reach the point $q$ by chronological products - trajectories of $\Sigma_M$ with piece-wise
constant controls; i.e., $q' \in {\Cal O}^-(q)$ if and only if $q\in {\Cal O^+(q')}$. Therefore, applying our Lemma~3 above we see that if some
general trajectory goes into the point $q$; i.e., there exists the trajectory $q(t), -\epsilon < t \leq 0, q(0)=q$ of $\Sigma_M$ containing two
points also connected by some general chronological product with negative arguments; then the following counterpart of the above theorem is true.

\proclaim{Theorem~2} If the control dynamical system $\Sigma_M$ satisfies the $LARC$ at the point $q$; i.e., there
exists the complete system - the set of generators -  $E_j,j=1,..,n', n'\geq n$ ($E_j$ generate $T_qM^n$) of
$T_qM^n$ consisting of some number of vectors $V_{u_i}, i=1,..., d'$, denoted by $W_i, i=1,..,d'$; and their
mutual commutators $W_i, i=d'+1,...,n'$
$$
W_i=[W_{\beta^i_k}, [W_{\beta^i_{k-1}}, [... ]]]; i=d'+1,...,n',
$$
then the set of points $q(t_1,...,t_N)$ (where all $t_i$ are negative) from which we can reach $q$ by general trajectories relative to this set of
generators is an open subset of $M^n$. If there exists a closed trajectory $q(t), 0\leq t \leq t_q$ going through the point $q$ then from some
open neighborhood $B^-_q$ of this trajectory we can reach $q$;
$$
q \in {\Cal O}^+(q') \quad \text{ for all } \quad q' \in B^-_q.
$$
\endproclaim

\medskip

\head 2. General trajectories and the involvement condition \endhead

Our definition of the general trajectory depends on the set of generators $V_{u_i}, i=1,...,d'$. The same
trajectory may be general relative to the one set of generators and not - relative to another set.

\proclaim{Definition} We say that the system $\Sigma_M$ satisfying the $LARC$ at some point $q$ is {\bf involved}
if for an arbitrary trajectory $q(t)$ going through this point there exists some set of generators $\{V_{u_i},
i=1,...,l'\}$ such that the trajectory $q(t)$ is general relative to this set.
\endproclaim

In general, the verification of the "involvement" condition is complicated and may require the solution of the
classification problem for vector fields on $M^n$.\footnote{Say, is the dependence $V_u$ on $u$ "stable" in some
sense or not?} Here we single out the case when the answer is easy.

Denote by $V(q)$ the set of all control vectors $V(q,U)=\{V_u, u\in U \}$ at this point, and by $L(q)$ the linear
subspace in $T_qM^n$ generated by $V(q)$. By $Con(q)$ we denote the convex hull of $V(q)$. We may assume that the
dimension ($d'$) of $L(q)$ is locally constant and the family of vector spaces $L(q)$ itself is smooth and locally
trivial - the union of all $L(q')$ over the points $q'$ from some small neighborhood $B(q)$ of $q$ - which we
denote by $L_{B}=\{ L(q')| q'\in B(q)\}$ - has a natural structure of the direct product $B(q)\times R^{d'}$ (is a
locally trivial vector bundle). The set $L_B$ is a subset of the restriction of the tangent bundle $TM$ over
$B(q)$, and admits natural coordinates $\{(x^1,...,x^n, y^1,...,y^n)\}$ in which the vector $y^j(\partial
/\partial x^j)(x)$ tangent to $M$ at the point with coordinates $x=(x^1,...,x^n)$ has coordinates
$(x^1,...,x^n,y^1,...,y^n)$.\footnote{For $d'\not= n$ the set $L_B$ is not contained in any coordinate subset of
the form $y^j=0$ in any coordinates.}

\proclaim{Definition} We say that the set of all control vectors $V_u, u\in {\Cal U}$ is {\bf open at $q$} if the
subset $\{ V(q')| q'\in B(q)\}$ is an open subset of $L_B$.\footnote{which is stronger than simply to require: "
$V(q')$ is an open subset of $L(q')$".}
\endproclaim

Let $V_i,i=1,...,N$ be some set of generators. Consider some family of linear transformations of $C(q'): R^N \to
R^N$ for $q'$ in a small neighborhood of $q$ given in the basis $\{V_i, i=1,...,N\}$ by a matrix
$C(q')=(c_i^k(q')); i,k=1,...,N$, and define the vector fields $\tilde V_k=c_k^i V_i$. Then, since $V(q')$ is
open, all vector fields $\tilde V_k(q')$ belong to $V(q')$ when $(c_i^k(q'))$ close to the identity matrix; i.e.,
they are control vectors of our control system $\Sigma_M$. In particular, if $C(q',\tau)$ is some continuous
family of linear transformations defined in some neighborhood of $q$ such that $C(q,0)\equiv Id$ then for small
$\tau$ vector fields $\tilde V_k(q',\tau) = c_i^k(q',\tau) V_k(q')$ belong to $V(q')$. We denote by $\tilde{\Cal
G}(\tau)(t_1,...,t_N)$ the corresponding chronological product defined with the help of these vector fields:
$$
\tilde{\Cal G}(\tau)(t_1,...,t_N) = P^{t_1}_{\tilde V_1(\tau)} \circ ... \circ P^{t_N}_{\tilde V_N(\tau)}. \tag 39
$$

Now we prove that it is possible to move the set of generators in an arbitrary way inside open $V(q)$, while
keeping the equation (29) of the zero set stationary. In order to do so, we define derivatives $U_i(\tau)$ of
$\tilde V_i(\tau)$ in the next lemma; where for simplicity we consider the moment $\tau=0$. The result, of course,
is valid for all $\tau$ sufficiently small.

\proclaim{Lemma~4} For an arbitrary vectors $U_i, i=1,...,N$ in $L(q)$ there exists a family of linear
transformations $C(q',\tau)$ in some small neighborhood of $B(q)$ of $q$ such that
$$
\delta (\tilde V_i(\tau)/\delta\tau)_{| \tau=0}= U_i, \tag 40
$$
and for all $i$,$|\alpha_i|>1$ it holds
$$
{{\delta \tilde V_i^{\alpha_i}(\tau)(x^j)}\over {\delta\tau}}_{| \tau=0} = 0, \tag 41
$$
where $C(q',\tau)$ is given by the matrix $(c_i^k(q',\tau)), i,k=1,...,N$, and $\tilde V_i(q',\tau) =
c_i^k(q',\tau) V_k(q')$.
\endproclaim

\demo{Proof} Indeed, by a direct calculation we see that
$$
{{\delta \tilde V_i^{\alpha_i}(\tau)(x^j)}\over {\delta\tau}}_{| \tau=0} =\sum\limits_{s=0}^{\alpha-1} (\tilde V_i^{s}\circ U_i \circ \tilde
V_i^{\alpha_i-s-1}). \tag 42
$$
We find $c_i^k(\tau)$ for each $i$ separately as follows: for a given $i$ find a coordinate system
$(x_i^1,...,x_i^n)$ in a neighborhood of $q$ such that the particular vector field $V_i$ with a given $i$ coincide
with the first coordinate field $\partial / \partial x_i^1$ of this system. In such coordinates the equation (42)
takes the form:
$$
\sum\limits_{s=0}^{\alpha-1} {{\partial}^s \over {\partial (x^1)^s}} ({{\partial c_i^k}\over{\partial \tau}}_{| \tau=0}(v^j_k)) = 0, \tag 43
$$
which obviously has the solution $c_i^k(q',\tau)$, say:
$$
c_i^k(q',\tau) =\delta_i^k + \tau ({{\delta c_i^k}\over{\delta \tau}})_{| \tau=0}, \tag 44
$$
and such that the products $c_i^k v_k^j$ are constant on $x_i$ coordinates:
$$
({{\delta c_i^k}\over{\delta \tau}})_{| \tau=0} (v_k^j)(q')\equiv (U_i^k v_k^j)(q), \tag 45
$$
where $v_k^s(q')$ and $U_i^k$ are coordinates of $V_k(q')$ and $U_i(q)$ in the coordinate system
$(x_i^1,...,x_i^n)$. To complete the proof it is sufficient to repeat our arguments for all $i$, define functions
$c_i^k$ in different coordinate systems $(x_i^1,...,x_i^n)$, and then rewrite them all (i.e., change coordinates)
in the same previously given initial coordinate system $(x^1,...,x^n)$. Note, that (45) in particular defines all
jets of the functions $c_i^k$, and (43) holds for all $\alpha$ with $|\alpha|$ arbitrarily big.
\enddemo

The Lemma~4 defines the variation on $\tau$ of the set of generators $\tilde V_i(\tau)$ in such a way that the
zero set $\tilde \Omega$ given by (26) for the chronological product (39) is defined by the first members of the
approximation sum for $\tilde{\Cal G}(\tau)$: the action of the member $d_{i\alpha}D^1_\alpha$ in the asymptotic
sum for the differential $d\tilde{\Cal G}(\tau)(\partial/\partial t_i)$ on the coordinate function $x^j$ does not
depend on $\tau$ if $|\alpha|>1$ due to (41) above; we have:
$$
{{\partial}\over{\partial \tau}}  d\tilde{\Cal G}(\tau)(\partial/\partial t_i)(x^j)_{| \tau=0} =
{{\partial}\over{\partial \tau}} (\sum\limits_{\alpha} d_{i\alpha}(t_1,...,t_N) D^1_\alpha(x^j))_{| \tau=0} =
{{\partial}\over{\partial \tau}} (\sum\limits_{|\alpha|=1} d_{i\alpha}(t_1,...,t_N)D^1_\alpha(x^j))_{| \tau=0} =
U_i(x^j). \tag 46
$$

As we wrote before, the differential of the chronological product may be not a surjection at the point
$t=(t_1,...,t_N)\in R^N$ if all partial derivatives $d\tilde{\Cal G}(\tau)(\partial/\partial t_i)$ at this point
vanish on some smooth function on $M^n$, which, without loss of generality we may assume to be some coordinate
function $x^j$. In other words, the image of the differential is contained in the subspace of $T_{q(t)}M^n$, which
at this moment ($\tau=0$) we may assume to be the tangent space to some coordinate submanifold in $M^n$ given by
$\{x^j=const\}$. At the origin $t=0$ - the point $q$ - this image (of $d\tilde{\Cal G}(\tau)$) always contains
vectors $V_i=V_{u_i}, i=1,...,N$. Therefore, we may assume that our coordinate systems $(x^1,...,x^n)$ are such
that
$$
V_{i}(x^j)(q)=0, \tag 47
$$
or, since the derivatives $U_i$ of $\tilde V_i(\tau)$ at $q$ also belong to $V(q)$, that
$$
U_i(x^j)(q)=0. \tag 48
$$

Hence, if $C(q',\tau)$ is defined with the help of the Lemma~4 for all $0<\tau<\epsilon$ and depends continuously
on $\tau$ we will have a continuous family of the sets of generators $\tilde V_i(\tau)$ such that in the equation
(29) defining the zero set $\tilde\Omega(\tau)$ of the chronological product $\tilde{\Cal G}(\tau)$ has the same
coefficients
$$
D^1_\alpha(\tilde V_1^{\alpha_1}(\tau)\circ ... \circ \tilde V_N^{\alpha_N}(\tau))\equiv const \tag 49
$$
for all $\alpha$. Therefore, $\tilde\Omega(\tau)$ does not depend on $\tau$ and is given in coordinates
$(t_1,...,t_N)$ by the same equation (29). Because of this for an arbitrary point $t^*$ from $\tilde\Omega(\tau)$
we may find a curve $t^*(\tau), 0\leq\tau<\epsilon$ issuing from this point $t^*(0)=t^*$ and outside
$\tilde\Omega(\tau)$ into the open domain $R^N$ of points with positive coordinates:
$$
t^*_i(\tau)>0 \quad \text{ and } \quad t^*(\tau)\not\in \tilde\Omega(\tau)\quad \text{ for }\quad \tau>0. \tag 50
$$
Next we define derivatives $U_i(\tau)$ by the following equation:
$$
t^*_i(\tau) U_i(\tau)(x^j)= - {{\delta}\over{\delta \tau}}t^*_i(\tau) d\tilde{\Cal
G}(\tau)(t^*(\tau))({{\partial}\over{\partial t_i}})(x^j), \tag 51
$$
which implies by direct calculations with the help of (49):
$$
{{\delta}\over{\delta \tau}} \tilde{\Cal G}(\tau)(t^*(\tau)) \equiv 0, \tag 52
$$
or by continuity that
$$
\tilde{\Cal G}(\tau)(t^*(\tau)) \equiv q(t^*). \tag 53
$$
We see that the continuous family of linear transformations $C(q',\tau), C(q',0)\equiv Id$ we just constructed,
move the vector fields $V_i(q')$ - the set of generators we have - into another set $\tilde
V_i(q',\tau)=c_i^k(q',\tau) V_k(q')$ of vector fields such that:

1) when $V(q)$ is open, then for $\tau$ small enough the vector fields $\tilde V_k(q',\tau)$ belong to our control
set $\{V_u(q')| u\in {\Cal U}\}$,

2) the point $q^*=q(t^*)$ which is the image under the chronological product ${\Cal G}(t^*)$ of some point from
the zero set $\Omega$ now equals the chronological product $\tilde{\Cal G}(\tau)(t^*(\tau))$ relative to another
set of generators, and is the image of the point $t^*(\tau)\in R^N$ which does not belong to the zero set
$\tilde\Omega(\tau)$ of this later chronological product; i.e., the differential of $d\tilde{\Cal
G}(\tau)(t^*(\tau))$ is the surjection on some small neighborhood of $q^*$.

In particular, all coordinate subspaces in $R^N$ defined by some number of equalities of the type $t_i=0$ belong to the zero set $\Omega$. For
instance, an arbitrary trajectory $q(t), 0\leq t$ of $\Sigma_M$ which is the solution of $\dot{q}(t)=V_{u(t)}(q(t))$, is the image of the first
coordinate line $\{t_i=0, i>1\}$ in $R^N$ under the chronological product ${\Cal G}(t_1,...,t_N)=P^{t_1}_{V_1}\circ ...\circ P^{t_N}_{V_N}$ - if
only we may consider the vector field $V_{u(t)}$ as the first vector field $V_1$ in the construction of our chronological products. If so, then
the arguments above implies that $\Sigma_M$ is involved. We come also to the same conclusion if trajectories with arbitrary controls $u(t)$ have
"nice" (e.g., smooth) dependence on $u$ and may be approximated by trajectories with piece-wise controls; i.e., by chronological products. To
treat the general situation we introduce one additional property.

\proclaim{Definition} We say that $V(q)$ is {\bf ample} if it is open and the convex hull of $V(q)$ coincide with
$L(q)$. \endproclaim

Now take an arbitrary trajectory $q(t),0\leq t\leq\epsilon$ under control $u(t)$ issuing from the point $q$ in
which $\Sigma_M$ satisfies $LARC$. Since $LARC$ implies local transitivity for $\epsilon$ sufficiently small
$q(\epsilon)$ equals some chronological product ${\Cal G}(t_1,...,t_N)$. Hence, the point $q$ is connected with
$q(\epsilon)$ by some trajectory $\bar{q}(t)$ of $\Sigma_M$ - possibly different from $q(t)$ - with piece-wise
control function. The part of such trajectory corresponding to the first non-zero $t_i$ is the solution of
$\Sigma_M$ under some constant control: $\dot{\bar{q}}=\pm V_i(\bar{q}(t)), 0\leq t\leq |t_i|$. If $t_i$ is
positive, then this part of the trajectory $\bar q(t)$ belongs to the positive orbit of $q$. If $t_i$ is negative,
and $V(q)$ is ample, then $V_i$ is the linear combination of some other control vectors $\tilde V_k$ with positive
coefficients
$$
V_i(q)=\lambda^k \tilde V_k(q), \tag 54
$$
and using the $LARC$ at $q$ it is not difficult to find a point $q^*$ on the trajectory $\bar q(t)$ which equals the chronological product $\tilde
{\Cal G}(\tilde t_1,...,\tilde t_N)$ defined with the help of the new vectors $\tilde V_k$ and such that the first non-zero argument $\tilde t_k$
is positive. In both cases we obtain the trajectory going from $q$ to $q(\epsilon)$ such that some part of it is some coordinate interval with
positive coordinates. Now the Lemma~4 claims that $q(t)$ is the general trajectory relative to some set of generators.

In other words, we just proved the following result.

\medskip

\proclaim{Lemma~5} If $V(q)$ is ample, then the control system $\Sigma_M$ satisfying $LARC$ at $q$ is involved.
\endproclaim

\medskip

Combining this with the claim of the Theorems~1 and~2 we conclude the following.

\medskip

\proclaim{Theorem~3} If for the system $\Sigma$ the set of control vectors $V(q)$ is ample at some point $q$,
where the $LARC$ is satisfied, then every closed orbit of $\Sigma_M$ going through the point $q$ has some open
neighborhood $B^+_q$ contained in the positive orbit ${\Cal O}^+(q)$ of $q$, and another open neighborhood $B^-_q$
such that for every $q'\in B^-_q$ the point $q$ is reachable from $q'$. Their intersection $B_q$ is an open
neighborhood of the closed trajectory on which the system $\Sigma_M$ is controllable.
\endproclaim


\medskip

\head 3. The Closed Orbit Controllability Criterium \endhead

Now we are ready to prove our main result:

\medskip

\proclaim{Theorem~A} Let $\Sigma_M$ be some control dynamical system on the compact connected manifold $M^n$. If through an arbitrary point $p$ of
$M^n$ goes some closed orbit $p(t),0\leq t\leq t_p$ which is non-trivial (i.e., $t_p>0$) and contains some point $q(p)$ where $V(q(p))$ is ample
and the $LARC$ is satisfied, than the system $\Sigma_M$ is controllable on $M^n$.
\endproclaim
\demo{Proof} The proof is obvious: every such orbit has an open neighborhood $B_q(p)$ on which the system $\Sigma_M$ is controllable. The union of
all $B_q(p)$ covers $M^n$. Since $M^n$ is compact there exists some finite sub-covering
$$
M^n=\cup_{i=1,...,c} B_{q_i}. \tag 55
$$
Because $M^n$ is connected for two arbitrary point $q$ and $p$ we may find a finite sequence $B_{q_{i_1}},..., B_{q_{i_k}}$ from the finite
covering (55) such that $q\in B_{q_{i_1}}$, $p\in B_{q_{i_k}}$ and $B_{q_{i_s}} \cap B_{q_{i_{s+1}}} \not= \emptyset$. Take $q_0=q$; some $q_s \in
B_{q_{i_s}} \cap B_{q_{i_{s+1}}}$ and $q_{k+1}=p$. Since $\Sigma_M$ is controllable on each $B_{q_{i_s}}$ there exists a trajectory from $q_s$ to
$q_{s+1}$ for all $s=0,..., k$. Taking the composition of these trajectories we obtain the trajectory of $\Sigma_M$ going from $q$ to $p$. Since
$q,p$ where arbitrary this proves that $\Sigma_M$ is controllable on $M^n$.
\enddemo

\medskip

\proclaim{\bf The Closed Orbit Controllability Criterium} Note that if $\Sigma_M$ is controllable on $M^n$ then for an arbitrary point $q$ there
exists some closed trajectory of $\Sigma_M$ going through this point. Indeed, since $\Sigma_M$ is controllable there exists some its trajectory
$\gamma$ from $q$ to $q(-\epsilon), 0<\epsilon$ for an arbitrary trajectory $q(t), -\infty<t<\infty$ passing through $q$. The composition of
$\gamma$ and the interval $q(t),-\epsilon\leq t \leq 0$ gives us the closed non-trivial trajectory of $\Sigma_M$ through $q$. Thus, the Theorem~A
above provides us with the necessary and sufficient controllability condition for the control systems $\Sigma_M$ on compact manifolds with ample
set of control vectors and satisfying the $LARC$, which we call the closed orbit controllability criterium.
\endproclaim

\medskip

The Theorem~A works especially well when applied to control dynamical systems on surfaces, when the trajectory curve is also a hyper-surface;
i.e., locally divides the manifold. If, for instance, $M^2$ is a surface and $\Sigma_M$ has the finite set of generator vector fields $V_i$ such
that the corresponding dynamical systems $\Sigma_i=\{\dot q(t)=V_i(q(t))\}$ are stable; i.e., have non-degenerated zeros, then our closed orbit
controllability criterium essentially provides the complete answer to the controllability question. This answer may be given with the help of some
diagrams (pictures) describing iterations on the space of orbits of $\Sigma_M$, which we split into some number of intervals, corresponding to the
cell representation of the phase portraits of $\Sigma_i$. We describe this in the forthcoming paper [CM], where we pay the particular attention to
the projections of bi-linear systems in $R^3$ to the two-dimensional unit sphere $S^2$ of directions in $R^3$. As the example the next theorem
gives the sufficient controllability condition for such control systems: let
$$
\dot x(t)= A(u)x(t)= (A + u^1 B_1 + ... +u^d B_d)x(t), \qquad \qquad (\Sigma)
$$
be the control systems of linear equations in $R^3$, where the right-hand side - the linear operator $A(u)$ - is a
linear combination of some constant linear operators $A$ and $B_k$ which do not depend on the control parameter
$u\in R^d$. The bi-linear system $\Sigma$ defines the control system $\Sigma^{pr}$ on the unit sphere $S^2$ in
$R^3$, called the projected system, when we consider the evolution of unit directions of solutions $x(t)$ of
$\Sigma$. For each control parameter $u$ being fixed the phase portrait of the corresponding projected system
admits the natural cell decomposition into, so called Jordan cells. Then the application of our closed orbit
controllability criterium provides, among other, the following.

\medskip

\proclaim{Theorem~B} If for two control parameters $u$ and $v$ the right-hand side linear operators $A(u)$ and
$A(v)$ have complex eigenvalues $\lambda_C(u)$ and $\lambda_C(v)$ correspondingly, then the system $\Sigma$
satisfying $LARC$ is controllable if for the real eigenvalues $\lambda_R(u)$ and $\lambda_R(v)$ of $A(u)$ and
$A(v)$ it holds
$$
( \lambda_R(u) - Re(\lambda_C(u)) ) \quad ( \lambda_R(v) - Re(\lambda_C(v)) )  < 0.
$$
\endproclaim

\medskip

We conclude with the example of the control system on the two-dimensional torus which satisfies $LARC$ everywhere and satisfies our closed orbit
criterium, but is not controllable since the set of control vectors is not "involved".

\medskip

\proclaim{Example~2} Let $\{S^1(\phi), 0\leq \phi \leq 2\pi\} \subset R^2$ be the unit circle
$S^1(\phi)=(cos(\phi),sin(\phi))$ in the two-dimensional $(x,y)$ plane $R^2$. Define two vector fields on this
circle: $V_1(\phi)= - sin(\phi)(sin(\phi),-cos(\phi))$ which equals the unit tangent to $S^1$ vector field
multiplied by the $y$-coordinate and has two zeros $q_1=(-1,0)$ and $q_2=(1,0)$; and $V_2$ - any smooth vector
field vanishing outside  $1/2$-neighborhoods of $q_1$ and $q_2$, but which is not zero at these points. Since
these vector fields do not vanish simultaneously and $S^1$ is one-dimensional the control system $q(t)=V_i(q),
i=1,2$ satisfies $LARC$. On the other hand it is not controllable because in two point $p_1=(0,1)$ and
$p_2=(0,-1)$ set of all possible directions of trajectories going through these points is the unique vector
$(-1,0)$ pointing to the same direction ("to the left"); i.e., any trajectory starting from $p_1$ or $p_2$ stays
in the "left" part of the circle where $x<0$ and can not reach any point with positive $x$-coordinate: ${\Cal
G}^+(p_i)\subset \{x\leq 0\}\not= S^1$.

Next we multiply $S^1(\phi)$ by another circle $S^1(\psi)$, and denote by $V_3$ the unit vector field tangent to this $S^1(\psi)$. As the result
we obtain the two-dimensional torus $T^2(\phi,\psi)=S^1(\phi)\times S^1(\psi)$ and the control system on it defined by three vector fields
$V_i,i=1,2,3$. Easy to see that this system 1) satisfies $LARC$ everywhere, 2) for an arbitrary point $q$ there exists non-degenerated closed
trajectory of the system going through this point - the second factor $S^1(\psi)$, 3) the system is not controllable.
\endproclaim

\enddocument

\medskip


\Refs
\widestnumber \key {AAAAAA}

\medskip




\ref \key  AS \by A.~A.~Agrachev, Yu.~L.~Sachkov \book Control Theory from the Geometric Viewpoint \yr 2003
\endref

\ref \key ACK \by V.~Ayala, E.~Cruz, W.~Kliemann \paper Controllability of Bilinear Control Systems on the
Projective Space\jour \vol \yr 2009 \endref

\ref \key ASM \by V.~Ayala, L.~A.~B.~San Martin \paper Controllability of two-dimensional bilinear systems:
restricted controls and discrete-time \jour \vol \yr \endref

\ref \key BGRSM \by C.~J.~Braga B., J.~Goncalves F., O.~do Rocio, L.~A.~B.~San Martin \paper Controllability of
two-dimensional bilinear systems \jour \vol \yr \endref

\ref \key CK \by F.~Colonius, W.~Kliemann \book The Dynamics of Control \yr 2000 \ed Birkhauser \endref

\ref \key CM \by A.~Choque,  V.~Marenitch \paper Direction control of bilinear systems, in preparation. \endref



\ref \key E \by D.~L.~Elliott \book Bilinear Control Systems \yr 2007 \endref


\endRefs
\bye